\documentclass[a4paper,10pt]{article}
\usepackage{multirow}
\usepackage[dvips]{color}
\usepackage{amsmath}
\usepackage{amssymb}
\usepackage{textcomp}
\usepackage{array}
\usepackage{rotating}
\usepackage{vmargin} 
\setmarginsrb{3cm}{4cm}{3cm}{4cm}{0cm}{0cm}{0cm}{0cm}

\usepackage{graphics}

\title{On green routing and scheduling problem}
\author{Nora Touati-Moungla \\ 
LIX, \'Ecole polytechnique, 91128 Palaiseau Cedex, France. \\
Vincent Jost \\
LIX, \'Ecole polytechnique, 91128 Palaiseau Cedex, France.
}

\begin{document}

\maketitle

\begin{abstract}
The vehicle routing and scheduling problem has been studied with much interest within the last four decades. In this paper, some of the existing literature dealing with routing and scheduling problems with environmental issues is reviewed, and a description is provided of the problems that have been investigated and how they are treated using combinatorial optimization tools.
\end{abstract}


\section{Motivation}

During the last few years, Operations Research (OR) has extended its scope to include environmental applications \cite{BLO95, DAN97, SBI07}. Because in Europe 73\% of the oil is used for transportation purposes, the need to design efficient plans for sustainable transportation is evident. Advances in the transportation planning process and in the efficiency of transportation systems are key components of the development of sustainable transportation.

The routing of vehicles represents an important component of many distribution and transportation systems and has been intensively studied in the OR literature \cite{TOTB02}. In this paper, particular consideration is given to routing and scheduling models that relate to environmental issues, we will denote this class of problems as Green Routing and Scheduling Problems (GRSP). In \cite{SBI07}, the authors discuss different problems that relate to sustainable logistics, they focus on reverse logistics, waste management, and vehicle routing and scheduling problems. Some variants of routing and scheduling problems in connection with environmental considerations were described: the arc routing problem, which is considered as a major component in waste management, and the time-dependent vehicle routing problem which allows one to indirectly decrease gas emissions involved by transportation activity by avoiding congested routes.

We present in this paper some general tools for transportation decision-making under assumptions related to economic, environmental and social considerations. An exhaustive review of sustainable transportation problems and their treatment by OR tools is out of scope here given the generality of this area of research. Most of the work in this domain is still very much in development and some applications have only just started. The aim of this survey is to provide a clear presentation on how combinatorial optimization can contribute to sustainable transportation as well as a comprehensive survey covering all known green routing and scheduling problems and their variants. Therefore, we list some GRSPs that are studied in the literature \cite{SBI07} and identify some other problems which can be added to this class, we describe some models and solution methods that can be exploited for these problems, and we expose some multiobjective optimization methods which are essential for solving these particular problems.

\subsection{Green routing and scheduling} 
\label{S_GVRS}

Usually, routing and scheduling models are concerned with objectives of minimizing economic costs, but due to growing concerns about public health, global warming and economic safety, it is necessary to consider in the cost function the factor of environmental and social costs. The additional environmental and social constraints and objectives that must be taken into consideration often make the problem more difficult to both model and solve. We study in this paper the best known routing and scheduling related problems arisen from sustainable transportation field: 

\begin{itemize}
 \item \textit{Routing of Hazardous Materials (RHM):} The objective is the minimization of the risk on the population and the environment caused by the transportation of hazardous materials. The research area that investigates this problem is most advanced, while less research taking into account explicitly environmental impacts is dedicated to the other problems cited bellow. 
 \item \textit{Routing and Scheduling in Time-Dependent Environment (RS\_TDE):} This class of problems contributes indirectly to reduce vehicles gas emission. The main objective is the minimization of a more realistic travel time by avoiding congested routes.
 \item \textit{Waste Collection Vehicle Routing Problem (WCVRP):} This problem is a major component of waste management.
 \item \textit{Multi-Modal Vehicle Routing problem (MMVRP):} This problem permits to manage many transportation modes and allows to perform priority-based routing for clean transportation (as rail transport).
 \item \textit{Dial-a-Ride Problem (DARP):} This problem contributes indirectly to decrease the global taxis gas emission by promoting grouped transportation (grouped taxis) for decreasing the transportation fleet size and routes congestion, particularly in large cities.
 \item \textit{Pick-up and Delivery Vehicle Routing Problem (PDVRP):} Permits the integration of the backward flow of waste in distribution systems for recycling for example.
 \item \textit{Energy Routing Problems (ERP):} This problem is one of the least studied one in the context of sustainable transportation. It permits to promote the use of electric vehicles by maximizing the vehicle autonomy.
 \item \textit{Air Traffic Control (ATC):} This problem can contribute to decrease the fuel consumption in planes.
\end{itemize}

These problems are very different from their structure, their contribution to transportation sustainability and their dedicated models and solution methods. Works on these problems are often unbalanced, this generally depends on the problem characteristics, for example RHM includes many aspects as risk definition and model, risk equity and uncertainty while the major aspect of RS\_TDE is the travel time.

\subsection{Limitations} 

A very large literature cover routing and scheduling problems presented in this survey, we do no give mathematical models of the presented problems since they are already well defined in the literature. However, an exception concerns energy routing problems for which no model is given in the literature. For each problem, we describe its characteristics, how it contributes to sustainable transportation, some models and resolution methods related to problems taking into account environmental considerations only. For example, in pick-up and delivery vehicle routing problem, the pick-up can concern materials or goods, we only interest to works in the literature which consider the pick-up of waste for recycling for example.

\subsection{Structure of the survey}

This paper is structured as follows. Section \ref{S_VRP} describes the classical vehicle routing and scheduling problem. In Section \ref{S_GVRS}, we describe the problems cited above, why they can be considered as green routing and scheduling problems, their particularity compared to the standard vehicle routing and scheduling problem, some related models and solution methods. We address in Section \ref{S_GVRSCP} the characteristics of green routing and scheduling problems and some classification schemes. Finally, we conclude in Section \ref{C_C}.

\section{A tour d'horizon of vehicle routing and scheduling models}

\subsection{Vehicle routing and scheduling problem}
\label{S_VRP}

The Vehicle Routing Problem (VRP) can be stated as follows: Consider a fleet of $K$ identical vehicles of fixed capacity $C$ available at a given depot to serve a set of costumers with fixed demand. We are given an oriented graph $G = (N,A)$, where $N$ is the set of nodes including the costumers and the depot, and $A$ the set of arcs connecting the nodes.  Each arc $(i,j)$ is associated with a cost $c_{ij}$ and each costumer $i \in N$ has a demand $d_i$. The goal is to find a set of minimum cost vehicle routes that service every costumer such that:
\begin{itemize}
 \item Each vehicle route originates and terminates at the depot.
 \item Each vehicle services one route and each costumer is visited by exactly one vehicle.
 \item The demand of each costumer is satisfied and the capacity of each vehicle is not exceeded.
\end{itemize}

The VRP is an important sub-problem in a wide range of distribution systems and a lot of effort has been devoted to research on different variants of this problem. Indeed, in practice, additional constraints or changes in the structure of the basic model are taken into account. We cite for example \textit{the VRP with time windows} which involves time window constraints restricting the times at which a customer is available to receive a delivery, in \textit{the VRP with Pick-up and Delivery}, each vehicle must visits the pick-up location before the corresponding delivery location, in \textit{the VRP with Backhauls} customers can demand or return some commodities, in \textit{the multiple depot VRP} the company may have several depots from which it can serve its customers and in \textit{the open VRP} each vehicle is not required to return to the central depot after visiting the final customer (see \cite{TOTB02} for a description of these VRPs). Many other variants of the VRP exist, however we focus in this work on vehicle routing and scheduling problems taking into account environmental considerations.

\subsection{Toward a general classification}
\label{S_TGC}

In this section we provide a classification of the models for some extensions of the VRP.
Although we do not want to provide a general nor a detailed classification
of all VRP models, our classification might inspire further research aiming at doing so.

Our motivation for proposing such a classification is to provide a
framework for describing and studying green variants and extensions of the VRP. 
This will be helpful for our exposition and might also be helpful for future works on 
(green) VRP, to distinguish between new applications and new refinements of models.

The main motivation for this classification is to show that the model refinements
that arose related to some applications might often be useful in other contexts. 
Hence a customizable general framework could be helpful in studying which refinements/extensions would be relevant in analyzing a given practical context.

This classification is inspired form and resembles the celebrated $\alpha|\beta|\gamma$-classification of machine scheduling problems. This could be helpful in further attempts to provided a detailed, formal classification of (green) VRP problems.

\subsubsection{Description of the demands}
\label{SS_DD}

\begin{itemize}
\item \textit{Physical nature of demands} Do we transport people, goods, hazardous materials, waste.
\item \textit{Volumes necessary for the demands} The objects to be transported require some space.
\item \textit{On-line/off-line demands} The demands can be known in advance or progressively introduced.
\item \textit{Origins-destinations} The demands can concern the delivery of products to certain costumers (distribution systems), the transportation of products from an origin location to a destination location (dangerous materials), or the pick-up of goods or people from a given location to another one.
\item \textit{Node/Arc demands} Demands can be associated with nodes, arcs or both of them.
\item \textit{Incompatibilities of demands} Do some demands share the same vehicle ? if not, can they be transported independently using the same vehicle ? 
\end{itemize}

\subsubsection{Description of the transport means:}
\begin{itemize}

\item \textit{Transportation network(s)} Do we have the possibility to ship over roads, rails, rivers and seas, sky, or a combination of these ? 

\item \textit{Fleet(s):} Are we dealing with a uniform fleet or several types of vehicles ? For each type of vehicle, what are its intrinsic characteristics: Load capacity (volume, weight, type of shippable products), maximum speed, energy consumption and autonomy, polluting emissions (per kilometer or as a function of speed). 

\item \textit{Availabilities, fares and intrinsic congestion:} For vehicles and arcs (paths) of the network, what are the periods of availability, the prices, the durations ? Do these durations depend on the dates at which we plan to 

\item \textit{Conflicts and logistical congestion:} In case we are in charge of coordinating
planes, security norms forbid the routes to be too close to each other (in the time-space diagram). More generally, when we are in charge of ``large'' (relatively to the network) fleets (of any transport means), our decisions may impact the congestion of the network, hence deteriorating our own performances.   
\end{itemize}

\subsubsection{Description of the objectives}

Beyond classical economic indicators, transportation industry also needs to be evaluated
on an environmental basis. Some impacts should be evaluated on a worldwide basis, but some impacts are mainly local, so that the evaluation should be performed on each ecosystem/human district.

\begin{itemize}

\item \textit{Classical/economical objectives:} Total revenue, total distance, Number of serviced demands, number of vehicles used, total operating costs etc.

\item \textit{Climate/Worldwide and sustainability objectives:} Total emissions of Green-House gazes, Total use of each primary/renewable source of energy.

\item \textit{Environmental(ecosystemic)/regional fairness and health objectives:} Local emissions in pollutants, impacts on human and animal populations, quality of air, water and soils, risks related to accidents.

\end{itemize}

\section{Handling the variety of demands}

As discussed on section \ref{SS_DD}, the characteristics of demand lead to a much different routing problems. Demands can concern a transportation request from a given location to another one (transportation of persons, transportation of hazardous materials), where routes are defined as shortest paths. When demands concern distribution requests where a set of costumers have to be serviced by a vehicle, we deal with vehicle routing. Demands can also be defined as pick-up and delivery system in which goods/people must first be picked up at a specific location and then be delivered elsewhere.

Time windows are generally imposed to restrict the time of the start and end of service or to restrict the time during which some road segments (arcs in the network) can be used, this last restriction is generally observed in transportation of hazardous materials.

Classically, deterministic demands are considered in the literature. Demands are known in advance and the routes are computed before the system starts to operate. In practice, demands are mainly stochastic since new demands can arrive during the service time and the routes have to be updated on-line (transportation of persons).

We present in the next sections how an adapted management of demands can lead to significants environmental savings.

\subsection{Dial-a-ride problem}
\label{SS_DARP}

The Dial-a-Ride Problem (DARP) consists of designing vehicle routes and schedules for $n$ costumers from their pickup point to their delivery point. The costumer requests the service by calling a central unit and specifying the origin and destination points, the number of passengers and some limitations in service time (the earliest departure time for example). The transport is supplied by a fleet of $m$ identical vehicles based at the same depot. The aim is to plan a set of minimum cost vehicles routes capable of accommodating as many requests as possible, under a set of constraints. The DARP can be static or dynamic. In the first case, the costumer asks for service in advance and the vehicles are routed before the system starts to operate. In the second case, requests are gradually revealed throughout the day and vehicle routes are adjusted in real-time to meet demand.

\subsubsection{Environmental contribution}

The main original target of the DARP is to offer the comfort and flexibility of private cars and taxis at a lower cost. This problem is suited to service sparsely populated areas, to low demand periods or to special classes of passengers with specific requirements (elderly, disabled). In addition, this problem considers indirectly environmental savings, indeed in opposition to individual taxis, the grouped ones can decrease the traffic density particularly in large cities. 

\subsubsection{Related works}

The DARP is characterized by multiple objectives such as the maximization of the number of costumers served, the minimization of the number of vehicles used and the maximization of the level of service provided on average to the costumer (costumer waiting time, total time spent in vehicles, difference between actual and desired drop-off times). The DARP can be formulated as multiobjective mixed integer program. Exact algorithms for the single-vehicle DARP have been developed in \cite{DES86, PSA80}. Recently, a branch-and-cut algorithm has been proposed in \cite{COR06a}. Heuristics and meta-heuristics are proposed for dealing with the dynamic problem with time-dependent network \cite{GIA03}. For a recent overview of the DARP, see \cite{COR06}.

In \cite{PAR09}, the authors propose a heuristic two-phase solution procedure for the dial-a-ride problem with two objectives. The first phase consists of an iterated variable neighborhood search-based heuristic, generating approximate weighted sum solutions and the second phase is a path relinking module, computing additional efficient solutions. 

\subsection{Pick-up and delivery vehicle routing problem}

In the Pick-up and Delivery Vehicle Routing Problem (PDVRP), a set of routes has to be constructed in order to satisfy transportation requests. A fleet of vehicles is available in a central depot to operate the routes and each vehicle has a given capacity. Each transportation request specifies the size of the load to be transported, the locations where it is picked up and the locations where it is delivered. Each load has to be transported by one vehicle from its set of origins to its set of destinations without any transshipment at other locations. The DARP (Section \ref{SS_DARP}) generalizes the PDVRP \cite{COR06}, the main difference between these problems is the human perspective; the level of service criteria is more important in the DARP.

\subsubsection{Environmental contribution} 

As mentioned in section \ref{S_CARP}, more and more countries have devoted considerable investments to waste reduction and material recycling. The existence of a backward flow of objects to be collected, stored, disassembled and recycled makes unprofitable to manage separately the forward flow of goods, from the producer to the  consumer, and the backward flow of waste or used-up devices, from the consumer to recycling or dumping facilities. In addition, when the reuse of products and materials becomes cheaper than simply disposing them, both of the opposite flows concern the producer, instead of being managed by independent subjects. When pick-ups concern waste, the PDVRP can be considered as a green routing and scheduling problem. This model derive from the development of reverse logistics, which consists of the efficient integration of the forward flow of goods with the backward flow of waste \cite{FLE97}. In their survey \cite{SBI07} underline the importance of reverse logistics in green Logistics, but the transportation aspect was not discussed. 

\subsubsection{Related works}

A comprehensive survey on the PDVRP can be found in \cite{PAR08} where different variants of the problem, models and resolution methods are presented. As our knowledge, no work in the literature treat real world pick-up and delivery problem in the context of recuperation of waste for recycling. Some applications of this problem can be (1) The door-to-door delivery of mineral water bottles and the simultaneous collection of empty bottles, (2) The laundry service for hotels (collecting dirty clothes and delivering clean clothes), and (3) Medical waste.

It is important to attach more interest to real problems for evaluate economical and environmental savings induced by these systems. The study of economical impacts of the integration of waste collection with products distribution can encourage industrials to recuperate the unused waste of their products and permits to the reduction of amounts of waste treated by municipalities and environmental saves.

\subsection{Waste collection vehicle routing problem}
\label{S_CARP}

Waste Collection Vehicle Routing Problem (WCVRP) can be classified as variation of the VRP but with additional constraints. The major difference between WCVRP and the classical VRP are landfills constraints. When a vehicle is full, it needs to go to the closest available disposal facility. Each vehicle can make multiple disposal trips per day. Three categories of waste are considered in the literature: commercial waste (involves servicing customers such as restaurants and small office buildings), residential waste (involves servicing private homes) and roll-on-roll-off waste (commonly used for construction site waste), these three categories bring about three different waste collection strategies. While works on the VRP consider the major objective of minimizing the travel cost, this problem also considers route compactness (a solution in which many routes cross over each other is less compact than one in which no routes overlap) and work balancing among vehicles.

\subsubsection{Environmental contribution}

In recent years waste management has become an area of concern for municipalities worldwide due to population growth, environmental concerns and the progressive increase in waste management cost. Waste collection is one of its main components. Note here that authors in \cite{SBI07} have discussed the importance attached to waste management and collection in terms of the ``green logistics'' agenda.

\subsubsection{Related works}

The particularity of residential routes compared to commercial and roll-on-roll-off ones is the mandatory adherence to driving on one side of the street. Unlike drivers on commercial or roll-on-roll-off routes, those on residential routes are permitted to serve only customers on the right side of the street. Very few exceptions are granted for alleys and one-way streets. Commercial and roll-on-roll-off waste collection differ principally from the size of the container. 
 
Routing problems in waste collection applications cannot typically be modeled with a unique routing model. As commercial and roll-on-roll-off waste routing consists of point-to-point collections, it can be modeled by node routing (VRP) models. Commercial waste routing problem can be characterized as a VRP with time windows (VRPTW) (Section \ref{S_VRP}) since commercial waste collection stops may have time windows. A special VRP variant known as rollon-rolloff vehicle routing problem was dedicated to roll-on-roll-off waste. However, residential routes require arc routing models, where costumers are located on the arcs. The periodic vehicle routing problem is also studied in the context of waste collection when collection operations are periodic on the time horizon. As these problems are computationally very hard, and can not be solved by optimal (exact) methods in practice, heuristics are used in this purpose. \\

\textbf{Vehicle routing problem with time windows:} The majority of papers in the literature are case study papers, focusing on results obtained when algorithms are applied to real-world data \cite{NUO06, SAH05, TUN00}. Only a few of these papers report computational experience with publicly available waste collection test instances \cite{KIM06}. More references are given in these cited papers.\\

\textbf{The periodic vehicle routing problem:} This problem has a horizon $T$, and there is a frequency for each customer stating how often within this $T$ period this customer must be visited. A solution to the problem consists of $T$ sets of routes that jointly satisfy the demand constraints and the frequency constraints \cite{ANG02}.\\

\textbf{The capacitated arc routing problem:} In this problem \cite{DRO00}, a fleet of vehicles, all of them located at a central depot and with a known capacity, must serve a set of streets network, with minimum total cost such that the load assigned to each vehicle does not exceed its capacity.\\ 

\textbf{The rollon-rolloff vehicle routing problem:} In this problem, tractors move large trailers between locations and a disposal facility. The trailers are so large that the tractor can only transport one trailer at a time \cite{BAL06, BOD00, DEM97}.

Waste collection real-life problems have almost been studied in a off-line context, where it is assumed that all data about the problem are known in advance. However, this is not necessarily the case when some information might not be readily available when the vehicles start their routes. When new information is available as the routes are executed, the problem becomes dynamic, this problem have not attracted yet the attention of waste collection research community.

\section{Handling congestion}

Most traffic networks, particularly road and air, face high utilization levels and congestion. As a result, the traffic conditions and its resulting variability can not be ignored in order to carry out a good quality route optimization. When taking into account congestion in the road context, RS\_TDE problems can be considered as major components for dealing with urban transportation problems associated with negative environmental impacts such as air pollution and noise. The problem of congestion in the airspace also cause a major concern given the forecast growth in aviation over the next decade. Air traffic delays due to congestion are a source of unnecessary cost for airline companies and passengers. Delays also have an environmental cost. Because of congestion, aircraft are often forced to fly far from the cruise altitude and/or the cruise speed for which they are designed, this results in unnecessary fuel burn and gas emissions.

\subsection{Routing and scheduling in a time-dependent environment}

The main difference between this class of problems and the classical VRP is the definition of travel time. When in classical VRP, the travel time is a function of the distance, in RS\_TDE problems, the travel time is variable and depends on many factors among which are weather conditions, congestion and the time of the day. We describe in this section three VRPs dealing with more realistic considerations of the travel time, the Time-dependent Vehicle Routing Problem (TDVRP), the Dynamic Vehicle Routing Problem (DVRP) and the Real-Time Vehicle Routing Problem (RTVRP). These problems are important not just because the consideration of the travel time variations affects considerably the objective values, but also because the best solutions known for non time-dependent problem are in general infeasible when applied in time-dependent world.

\subsubsection{Time-dependent vehicle routing problem}
\label{SS_TDVRSP}

When the VRP assumes that the costs or travel times are a scalar transformation of distance, the TDVRP is more adapted to real applications by taking into account variations of the travel time resulting from periodic cycles in the average traffic volumes. It is considered in this problem that the principal variation in travel time results from the time-of-day variation, the travel time between two points $i$ and $j$ is a function of the time of the day at the origin point $i$. A variety of models for the TDVRP are considered in the literature. We present briefly a classification of these models \cite{ICH03}:
\begin{enumerate}
 \item \textit{Basic Models (BM):} Time-dependency is integrated in the model using simple rules like multiplier factors associated with different periods of the day. Unfortunately, these assumptions are weak approximations of the real-world conditions where travel times are subject to more subtle variations over time.

\item \textit{Models based on Discrete Travel time and Cost Functions (MDTCF):} In this kind of formulations, the horizon of interest is discretized into small time intervals. The travel time and cost functions for each link are assumed to be step functions of the starting time at the origin node. However, the assumption that travel times vary in discrete jumps is just an approximation of real-world conditions since travel times change continuously over time. Many of these models are dedicated to time-dependent shortest path problem \cite{CHA97} and time-dependent traveling salesman problem \cite{MAL92}. In \cite{KOS93}, the authors consider the problem of path planning in networks including multiple time dependent costs on the links and use the dynamic programming algorithm principle.

\item \textit{Models based on Continuous Travel time and Cost Functions (MCTCF):} Continuous travel time functions seem to be more appropriate to model real-world conditions. Unfortunately, the models obtained are difficult to solve without simplifying assumptions, so these models consider again an approximation of the real travel time variations. 

\item \textit{Queueing Models (QM):} Here, the traffic congestion component is based on queueing theory. This allows one to capture the stochastic behavior of travel times by generating an analytical expression for the expected travel times \cite{VAN07}. 

\end{enumerate}

Works in this field show experimentally that the total travel times can be improved significantly by explicitly taking into account congestion during the optimization. Very few comparative framework on different models and solution methods are found. We show on Table \ref{TAB_TDVRP} some works on the TDVRP.

\begin{table}[!h]
\scriptsize
\centering
\caption{Some works on the time-dependent vehicle routing problem} 
\label{TAB_TDVRP}
\begin{tabular}{p{3cm}p{1cm}p{5cm}p{5cm}}
\\
\hline
References & Models & Models characteristics & Solution methods\\
\hline

Brown, Ellis, Graves and Ronen (1987) \cite{BOW87} & BM & A solution where travel time fluctuations are ignored is first produced. Then, the loads for each truck are resequenced ``manually'' to take into account various factors such as traffic congestion during rush hours, road and weather conditions. & A collection of integer programming methods.\\
\hline
Malandraki and Daskin (1992) \cite{MAL92} & MDTCF & The problem is formulated as a mixed integer programming problem. & Nearest-neighbour (greedy) heuristic is proposed, as well as a branch-and-cut algorithm for solving small problems with 10-25 nodes.\\
\hline
Hill and Benton (1992) \cite{HIL92} & MDTCF & The model was based on time-dependent travel speed. & Experimentations based on a small example with a single vehicle and five locations are given. \\
\hline
Ichoua, Gendreau and Potvin (2003) \cite{ICH03} & MDTCF & The model was based on time-dependent travel speed. & A taboo search heuristic is proposed and experimentations are performed on Solomon's 100-costumers problems. \\
\hline
Woensel, Kerbache, Peremans and Vandaele (2007) \cite{VAN07} & QM & - & Both the static and the dynamic TDVRP were solved using ant colony optimization.\\
\hline
Hashimoto, Yagiura and Ibaraki (2008) \cite{HAS08} & MDTCF & Travelling time and cost functions values are time-dependent. & A local search algorithm.\\
\hline
Donati, Montemanni, Casagrande, Rizzoli and Gambardella (2008) \cite{DON08} &MDTCF & The model was based on time-dependent travel speed. & A multi-Ant Colony System.\\
\hline
\end{tabular}
\end{table}

\subsubsection{Dynamic vehicle routing problem}

The DVRP is the dynamic counterpart of the VRP, where information relevant to the planning of the routes can change after the initial routes have been constructed. This class of problems have arisen thanks to recent advances in communication and information technologies that allow information to be obtained and processed in real time. 

Traditionally, the DVRP is solved in one of two ways: when problem parameters become (customers, speed, travel time) known throughout the run, oblivious online algorithms are used. Alternatively, when all parameters are available but with uncertainties in their properties, stochastic optimization is used, which build the routing plan a priori, and then modifies it when changes in parameters properties occur. Most research in this area has focused on dynamic routing and scheduling that considers the variation in customer demands. However, there has been limited research on routing and scheduling with congestion and travel times variation. We present below some models for the DVRP proposed in the litterature:

\begin{enumerate}

 \item \textit{Models based on Simulations (MS):} A vehicle routing and scheduling plan is obtained by revisiting the vehicle routing and scheduling plan computed on the previous day by using real-time information of present link travel times, whenever a vehicle arrives at a customer. This real-time information is provided by dynamic traffic simulation based on the current conditions of the day \cite{TAN04}. 

 \item \textit{Queueing Models (QM):} To capture travel times, these models introduce a traffic congestion component based on queueing theory. A major advantage of using these queueing models is that the real-life physical characteristics of the road network can be mapped immediately into the parameters of the queueing model. Moreover, the inherent stochasticity of travel times can explicitly be taken into account via the analytical queueing models \cite{VAN08}. 

 \item \textit{Stochastic Models (SM):} The travel times are subjected to stochastic variations \cite{POT06}. 
\end{enumerate}

Due to the complexity of this problem, heuristic methods are often used to obtain good solutions, tabu search algorithms were used to solve QM \cite{VAN08}, genetic algorithms to solve MS \cite{TAN04} and local search heuristics to solve SM \cite{POT06}. A discussion of the DVRP solution methods can be found in \cite{PSA95}.

\subsubsection{The real-time vehicle routing problem} 
 
A new generation of VRP are proposed in the literature for vehicle routing in more realistic settings. The RTVRP considers more accurate information about the travel times compared to the DVRP by considering real-time variations in travel times \cite{OKH09, ZHU00}. Dynamic vehicle routing models imply that a vehicle en route must first reach its current destination and only after that can it be diverted from its route. However, an unpredicted congestion or other traffic impediment can be encountered on the way to its immediate destination. Using mobile technology, vehicle routing can be modeled in more realistic settings:

\begin{itemize}
 \item Allows locating vehicles in real time.
 \item Enables the online communication between the drivers and the dispatching center.
 \item Capable to capture varying traffic conditions in real time and in the short run predict with high accuracy the travel time between a pair of nodes.
\end{itemize}

All these factors allow to send new instructions to drivers at any time, regardless of their location and status. These modeling approaches enable a better approximation of the real-world conditions.

\subsubsection{Synthesis}

It was established that the highest emissions of carbon dioxide occur in congested, slow moving traffic. RS\_TDE models can contribute indirectly to decrease fuel consumption and gas emissions \cite{SBI07}. In reality, the links have different combinations of congestion levels, and delays associated with road furniture such as traffic lights and roundabouts, and road topography and geometry such as inclines. This causes speed variations (resulting from acceleration and deceleration) and therefore produces different times over links with the same road category and distance. In addition, these speed variations would lead to fuel consumption and therefore gas emissions variations. Therefore, there is a need to calculate vehicle routes which minimize gas emissions, not just to calculate routes minimizing by time or distance. 

\subsection{Air traffic control}

Air traffic management consists of two important components: the traffic planning and the traffic control. Traffic planning deals with the balance between demand and the available capacity and traffic control has to guide aircraft safely to their destinations. 

An aircraft conflict occurs when the distance between two or more aircrafts falls below a given threshold. In this case, a minimum separation is required. Aircraft conflict detection and resolution has been widely studied in the literature \cite{KUC00}. In \cite{PAL02}, the authors propose an integer programming model that minimizes the maximum deviation in the changes made, by assuming that aircraft can perform either a speed change or a heading change. Authors in \cite{CLA08} studied the traffic control problem by maintaining separation while considering associated fuel costs with any heading deviation or speed changes. Safety requirements are considered as hard constraints that must be maintained. The objective function focuses on minimizing fuel costs, and hence the resulting environmental impact. In \cite{BER00}, the authors address the problem of determining how to reroute aircrafts in the air traffic control system when faced with dynamically changing weather conditions. The overall objective of this problem is the minimization of delay costs. 

The multiobjective aspect of the problem was considered in \cite{TIA10} where a multi-objective genetic algorithm has been designed to solve the model. Three objectives have been considered: the minimization of the sum of the flights which exceeds the capacities of all the sectors, the minimization of the sum of the maximum flights which exceeds the capacities of all the sectors and the minimization of the total delay time including the ground delay time and air delay time. In \cite{VIJ08} the authors take into account the weather-related flight delays and investigate the problem of generating optimal weather avoidance routes under hazardous weather conditions. The proposed model minimizes the fuel usage, weather conditions, customer comfort and traffic density and is solved using a hybrid ant colony optimization method with a multi-agent approach.

\section{Toward cleaner transportation means}

Clean transportation means are generally considered as those which include public transportation, railway freight, and cycles. They can be defined as modes using sustainable energy sources such as renewable energy. Nuclear energy for electric drive systems is often referred to as sustainable, but this is controversial politically due to concerns about peak uranium, radioactive waste disposal and the risks of disaster due to accident, terrorism, or natural disaster as the case on March 2011 at Japan.

\subsection{Routing of electric vehicles}
\label{SS_ERP}

Electric vehicles are generally considered as the cleanest transportation option, due to their zero local green house gas emissions and noise, particularly in large cities. The succeeding of a transition from a conventional gasoline based transportation system towards a sustainable way of transportation, depends on a quite number of critical factors. The substitution of conventional vehicles through electric and/or hybrid vehicles involves economical, environmental and social aspects. The purpose of EVM is to encourage the transition to EV use and to expedite the establishment of a convenient, cost-effective, EV infrastructure that such a transition necessitates. Whereas the development of EV through battery autonomy grow extensively, very little research is dedicated to EV routing, recharge stations localization, vehicles redistribution and energy management. These major aspects represent challenging optimization problems.

\subsection{Multi-modal routing problem}
\label{SS_MMRP}

This problem is defined as follows: given a set of origin-destination transport requests, one must optimally route these requests in a multi-modal network including a heterogeneous set of transportation services. These services are generally classified by according to two main characteristics: the departure time and the cost function. By the departure time characteristic, we can differentiate between timetable services (rail and short sea shipping) and time-flexible services (trucks). The cost (and duration) of routes depends on the departure time, the transportation mode, the distance and the waiting time in transshipment nodes. These constraints make the problem even more difficult. 

\subsubsection{Environmental contribution}

Negative impacts caused by the transportation activities such as gas emission and noise can be reduced by using cleaner and silent alternative transportation modes. Multi-modal transportation strategies are studied in the operations research literature through the multi-modal routing problem.

\subsubsection{Related works}

In \cite{CHA08}, the authors consider a multiobjective multimodal multicommodity flow problem with time windows and piecewise linear concave cost functions. Based on Lagrangian relaxation technique, the problem is broken into a set of smaller and easier subproblems and the subgradient optimization procedure is applied to solve the Lagrangian multipliers problem. Authors in \cite{MOC08} proposed an origin-destination integer multi-commodity flow formulation with non-convex piecewise linear costs and use column generation based heuristic that provides both lower bounds and good quality feasible solutions. The author in \cite{MIN91} deals with two objectives: the cost and the risk and develops a chance-constrained goal programming method to solve the problem.

Multi-modal transportation models need to define the optimization methods from the view of economic and environmental performances. When we are interested in multi-modal transportation, it is essential to make an adapted choice of transportation modes in such a way that the environmental impacts of the transportation system are minimized. In \cite{BON04}, the authors underlined the importance of choosing the transportation modes but little work is known on the calculation of transportation cost, taking environmental impacts into account. In \cite{ANC05}, the transportation mode is considered within the framework of a green supply chain, the environmental impacts of the means of transport are integrated into the model as a cost aspect.  

In multi-modal transportation, specific characteristics must be taken into account for each transportation mode. For example, the impact of variation in travel time on rail transport links and road transport links is very different, as rail transport mode is not subject to congestion. For determining efficiently the transportation mode for minimizing the environmental impacts, it is important to consider stochastic travel time. The multi-modal transportation network can be assumed in this case to be mixed, where the travel time on some arcs is stochastic.  

\section{Multi-objective assessment of sustainability and quality}

\subsection{Environmental impacts of transportation}

Transport activities have resulted in growing levels of motorization and congestion. As a result, the transportation sector is becoming increasingly linked to environmental problems. With a technology relying heavily on the combustion of hydrocarbons, notably with the internal combustion engine, the impacts of transportation over environmental systems has increased with motorization. This has reached a point where transportation activities are a dominant factor behind the emission of most pollutants and thus their impacts on the environment. 

Knowing how transport activities impact on the amount of CO2 and local pollutant emissions is an important step in the evaluation of system sustainability, but it is not always the final goal since the nature and the intensity of impacts on the population and the environment on local and global scale and on short and long time horizons are also key questions. Therefor, there is an urgent need of approaches and methodologies to determine the relationships between transport and the environment. Environmental transportation impacts, can fall within three categories:

\begin{itemize}
 \item \textbf{Direct impacts:} The immediate consequence of transport activities on the environment where the cause and effect relationship is generally clear and well understood.

 \item \textbf{Indirect impacts:} Include changes in land use and impacts on environmental resources such as habitat fragmentation on species viability over time or climate change. Indirect effects may have greater consequences than direct effects, but are not generally well understood.

 \item \textbf{Cumulative impacts:} Cumulative impacts take into account the varied effects of direct and indirect impacts on an ecosystem, which are often unpredicted.
\end{itemize}

Environmental impact assessments of transport is very complex and have to aggregate properly the whole range of impact varieties (residential households, surface water, flora and fauna, etc) and their relative importance by bringing together leading experts and practitioners on transport, air quality, climate change and related issues. The main environmental dimensions of transportation include \cite{OFF96}:

\begin{itemize}
 \item \textbf{Causes.} Represent factors that influence transport activities such as land use, demographics and economics. These mesures do not provide a great deal of information for estimating the environmental consequences of transportation, but they help to expalin the reasons why certain impacts may be increasing or decreasing.

 \item \textbf{Activities.} Involve a wide array of factors expressing the usage of transportation infrastructures and all the related services. Activities are related to \textit{infrastructure} (The environmental impacts of the construction and maintenance of transport infrastructure), \textit{vehicle manufacture} (The resources and energy consumed in the manufacturing process have an environmental impact), \textit{vehicle travel} (Environmental impacts the outcome of regular transport vehicle use, which varies by mode), \textit{vehicle maintenance} (The consumption and disposal of parts during maintenance can have environmental impacts) and \textit{vehicle disposal} (Once the useful life of a vehicle is over, its disposal may have some environmental impacts). Activities often have direct environmental consequences and the level of environmental damage associated with a specific activity varies bu location and over time.

 \item \textbf{Outputs.} Several factors are to be considered. The first outcome of transportation activities are gaz emissions. According to the geographical characteristics of the area where emissions are occurring (e.g. wind patterns) ambient levels are created. Once these levels are correlated with population proximity, a level of exposure to harmful pollutants can be calculated. 

 \item \textbf{Results.} They include all the health, environmental and welfare effects of the exposure to emissions from transportation activities and they are very difficult to measure.

\end{itemize}

The relationships between transport and the environment are also complicated by the different dimensions of impacts. Transport activities contribute at different geographical scales to environmental problems, ranging from local (noise and gaz emissions) to global (climate change). Currently, there is still a lack of coherent, consistent and applicable methodologies for assessing such impacts from transport interventions. In addition, there are limitations on available and consistent data to evaluate climate and air pollution impacts of transport and related investments.

\subsubsection{Local impacts}
Local concentrations of atmospheric pollutants in urban areas are an increasing problem for many countries.  Poor air quality causes premature mortality and respiratory diseases particularly in large metropolitan regions. Behind atmospheric pollutants, other local impacts result from transportation activity: 

\begin{itemize}

 \item Noise nuisances are a serious problem, particularly of road transport in densely populated areas. Noise associated with shipping has the potential to cause disturbance to marine animals. In aviation, standards and policy positions regarding aircraft noise have long existed and the issue of aircraft noise has been a driving issue in the location and operational parameters of airports. 

 \item The development of roads and highways is associated with a number of major environmental issues such as those associated with loss of land resources through construction of infrastructure.

 \item Wildlife habitat loss: Transportation activities create a barrier to wildlife movement causing the reduction of biological diversity, and the increasing threat of extinction \cite{ALE00}. As our knowlage there no work dealing with the impacts of transportation traffic in terms of gas emission on wildlife corridor.

\end{itemize}

Few quantitative impact assessment of transport policies have been published and there is a lack of a common methodology for such assessments. In \cite{SCH09}, the authors evaluate the usability of existing health impact assessments methodology to quantify health effects of transport policies at the local level using two simulated transport interventions. The first is based on the speed limit reduction and the second on the traffic re-allocation with an aim of improving the neighboring populations quality of life. Due to the lack of data on population health, the effects of these interventions however appeared low (5\%). In \cite{IYE09}, the author                                                                             proposes an economic model for quantifying health impacts of public policies, particularly in transport, in terms of dollar value. For more details on quantitative health impact assessment of transport policies, see \cite{VEE05}.

\subsubsection{Global impacts}

Global pollution is also an increasing problem. Adapting to climate change is an evolutionary process which requires adoption of longer planning horizons, risk management, and adaptive responses for making transportation activity more sustainable. The strategic examination of national, regional, and local networks is an important step toward understanding the global risks posed by transportation activities. 

Air pollution is the most important source of environmental externalities for transportation. Although the nature of air pollutants is clearly identified, the scale and scope on how they influence the biosphere are subject to much controversy. Air pollution costs are probably the most extensive of all environmental externalities of transportation. As all externalities, costs are very difficult to evaluate because several consequences are not understood and a monetary value cannot be effectively attributed. In general, the costs of air pollution associated with transportation can be grouped within economic, social and environmental costs.

\subsubsection{Evaluation of transportation usage impacts}

Total costs incurred by transportation activities, notably environmental damage, are generally not fully assumed by the users. The lack of consideration of the real costs of transportation could explain several environmental problems. For instance, external costs account on average for more than 30\% of the estimated automobile costs. If environmental costs are not included in this appraisal, the mobility cost is partially assumed by the users (e.g. fuel, licensing, insurance, etc.). However, environmental impacts is a cost mostly assumed by the society. Therefore, there is an urgent need of approaches and methodologies to determine the relationship between emissions, air quality and health impacts as well as their monetary quantification.

At European level, it has been recognized that strategic environmental assessment should be developed as an integral part of the decision-making process for policies, plans and programs \cite{EC96, PAL07}. In the case of transportation, evaluation of the environmental impacts begins with an identification of possible impacts on various aspects of the natural and social environment. The difficulties in arriving at a generalized framework for environmental evaluation arise at two levels. On the one hand, such framework should be easy to be understood and to employ by decision makers, while meeting a series of efficiency criteria, such as generality, independence, reliability, reasonable data needs, etc. On the other hand the evaluation framework should be realistic and coherent, able to deal with challenging issues, such as uncertainty, real-time, network effects and modal changes. 

\subsection{Routing problem for hazardous materials}
\label{RHM}

The transportation of hazardous materials (hazmat from here on) has received much interest in recent years, this results from the increase in public awareness of the dangers of hazmats and the enormous amount of hazmats being transported. The main target of this problem is to select routes from a given origin $s$ to a given destination $t$ such that the risk for the surrounding population and the environment is minimized, without producing excessive economic costs. The study of hazmat transportation problems can be classified in four main subjects: risk analysis \cite{ERK05}, routing and scheduling \cite{LIS91}, facility location \cite{LIS91}, and treatment and disposal of waste \cite{NEM99}, we focus in this section on routing and scheduling.

\subsubsection{Environmental contribution}

This problem is naturally a GRSP, since it contributes to minimizing the risks of release accidents on the population and the environment in hazmat transportation activities.

\subsubsection{Related works}

The difference between hazmat transportation and other transportation problems is mainly the risk. The risk makes this problem more complicated by its assessment, the related data collection and the solution of the induced formulations. 

\subsubsection*{Problem characteristics}

We address in the following two major particularities of RHM problem: risk assessment and risk equity:\\   

\textbf{(a) Risk assessment:} Although the fact that the major target of RHM is the minimization of the risk, there is no universally accepted definition of risk (for a survey on risk assessment, see \cite{ERK98}). The risk caused by hazmat transportation depends on many factors, the most important of which are: the risk categories (explosion, toxicity, radioactivity, etc), the transportation mode, the affected agents (population, territorial infrastructures and natural elements), the meteorological conditions and the temporal factor. It is pointed out in \cite{LIS91} that the evaluation of risk in hazmat transportation generally consists of the evaluation of the probability of an undesirable event, the exposure level of the population and the environment, and the degree of the consequences (e.g., deaths, injured people, damages). In practice, these probabilities are difficult to obtain due to the lack of data and generally, the analysis is reduced to consider the risk as the expected damage or the population exposure.

As the risk is a part of the objective function, it is quantified with a path evaluation function \cite{ERK98}. This function is not additive since the probability of a release accident on a link depends on the probability of a release accident on the traveled links of the path. This important property leads to non-linear integer formulations which can not optimized using a classical shortest path algorithm. Generally approximations are needed by considering additive functions (Usually considering independent release accident probabilities on links) for obtaining tractable models.  \\

\textbf{(b) Risk equity:} When many vehicles have to be routed between the same origin-destination nodes, these vehicles are routed on the same path, hence the risk associated to regions surrounding this path could be high. In this case, one may wish distribute the risk in an equitable way over the population and the environment. 
The computation of routes with a fairly distributed risk consists in generating dissimilar origin-destination paths, i.e paths which relatively don't impact the same zones. Solution approaches ca be classified in tow sets, \textit{resolution-equity-based methods} and \textit{model-equity-based methods}. In resolution-equity-based methods, equity constraints are taken into account in the resolution process. These methods are based on a dissimilarity index which permits to indicate when two paths are considered as dissimilar. We present on Table \ref{TAB_EQ} some of these methods.

\begin{table}[!h]
\scriptsize
\centering
\caption{Resolution-equity-based methods} 
\label{TAB_EQ}
\begin{tabular}{p{2.5cm}p{7cm}p{4cm}}
\\
\hline
Method & Principle & dissimilarity index \\ 
\hline
Iterative Penalty Method \cite{JOH92} &Compute iteratively a shortest path and penalize its arcs by increasing their weights for discouraging the selection of the same arc set in the generated paths set in the next iteration. & \\
\hline
Gateway shortest-paths method \cite{LOM93} &Generate dissimilar paths by forcing at each time a new path to go through a different node (called the gateway node). & The absolute difference between areas under the paths (areas between paths and the abscissa axis). \\
\hline
Minimax method \cite{KUB97} & Select $k$ origin-destination shortest-paths and select among them iteratively a subset of Dissimilar Paths (DP) by means of an index that determines the inclusion or not of candidate paths in DP. & The length of common parts between the paths.\\
\hline
$p$-dispersion method \cite{AKG02}& Generate an initial set $U$ of paths and determining a maximal dissimilar subset $S$, i.e., the one with the maximum minimum dissimilarity among its paths. & The length of common parts or the common impact zones between the paths.\\
\hline
\end{tabular}
\end{table}

Model-equity-based methods consists of taking into account equity constraints in the model formulation. In \cite{GOP90a, GOP90b}, the authors propose an equity shortest path model that minimizes the total risk of travel, while the difference between the risks imposed on any two arbitrary zones does not exceed a given threshold $\epsilon$. In \cite{CAR08} was proposed a multi-commodity flow model for routing of hazmat, where each commodity is considered as one hazmat type. The objective function is formulated as the sum of the economical cost and the cost related to the consequences of an incident for each material. To deal with risk equity, the costs are defined as functions of the flow traversing the arcs, this imposes an increase of the arc's cost and risk when the number of vehicles transporting a given material increases on the arc. \\

\subsubsection*{Routing and scheduling}

Transportation of hazardous materials is a complex and seemingly intractable problem, principally because of the inherent trade-offs between social, environmental, ecological, and economic factors. This problem is multiobjective in nature, since risk minimization accompanies the cost minimization in the objective function. In addition, other pertinent objectives can be considered as the travel time for minimizing the exposure of the driver to risk. Therefore, a set of alternative (Pareto optimal) solutions have to be computed (see section \ref{S_MOO}). Solution methods for hazmat routing can be classified in two categories: \\

\textbf{Local hazmat routing:} Consists of a one origin-destination hazmat routing and aims at selecting routes between a given origin-destination nodes for a given hazmat, transport mode and vehicle type. We present in Table \ref{TAB2} some works on local routing. \\

\begin{table}[!h]
\scriptsize
\caption{Local hazmat routing} 
\label{TAB2}
\begin{tabular}{lll}

\\
\hline

Author &  Objective & Method \\ 
\hline

Shobrys (1981) \cite{SHO81}             & Min. ton-miles traveled       & Weighting method  \\
					& Min. population exposure-tons &                   \\ 
\hline

Robbins (1981) \cite{ROB81}  & Min. the total length of shipment       & Weighting method  \\
& Min. the size of the population brought &                   \\
                    & into contact with the shipment          &                   \\ 
\hline

Current, ReVelle and Cohon  & Min. the population affected  & Weighting method  \\
(1988) \cite{CUR88}	    & around the path		     &		         \\ 
 & Min. the length of the path    &                    \\ 
\hline

Abkowitz and Cheng  & Min. ton-miles travelled        & Weighting method  \\
(1988) \cite{ABK88} & Min. population exposure-tons  &   \\ 
\hline
Turnquist (1993) \cite{TUR93}   & Min. of incident rates related to       	& Stochastic dominance   \\
				& the release of hazardous material		&				  \\
				& Min. of the population exposed to the risk	&				  \\
				& Min. of the route length			&				  \\
\hline
Karkazis and Boffey (1995) \cite{KAR95} & Min. expected damage effects & Branch-and-bound   	\\
	   
\hline
I. Giannikos (1998) \cite{GIA98}& Min. Operating cost                     & Goal Programming,    \\
		    & Min. Total perceived risk               & Penalty Functions             \\
 		    & Min. Maximum individual perceived risk  &                      \\
 		    & Min. Maximum individual disutility      &                               \\ 
\hline
Zografos and Androutsopoulos      & Min. Total travel time         & Objectives aggregation,  \\
(2004) \cite{ZOG04}    & Min. Total transportation risk & Insertion heuristic   \\
\hline

\end{tabular}

\end{table}

\textbf{Global hazmat routing:} A substantial work in the literature focuses on the selection of a single commodity routes between only one origin-destination pair. In practice, a better suited model is global routing, where different hazmats have to be shipped simultaneously among different origin-destination pairs. 
In \cite{ZOG89}, a multiobjective routing model that considers equity constraints is proposed. The model considers the following criteria: the general risk, the risk of special population, the travel time, the property damage, and the risk equity which is imposed using capacity constraints on the network links. The obtained model is equivalent to the capacitated assignment problem and a goal programming method is used to solve it. A pertinent model for global routing was proposed in \cite{CAR08} and described in section \ref{RHM}-(b), unfortunately the multicommodity flow model is mono-objective in nature. As our knowledge, the multiobjective aspect of this model is not yet studied for the considered problem.

\subsection{Resolution methods}
\label{S_MOO}

The GVRSP was defined in Section \ref{S_GVRS} as a classical VRP with additional environmental and social objectives and constraints, this radically changes the problem structure, so different and generally dedicated solution methods are developed. As a consequence, we envision a new era in which optimization systems will not only allocate resources for optimizing economical costs: they will react and adapt to external events efficiently under environmental constraints and objectives. The GVRSP can be defined as the computation of a set of origin-destination paths optimizing a set of objectives $\mathcal{F}$ by satisfying a set of constraints $\mathcal{C}$. We summarize in Table \ref{TABLE3} some general characteristics of the presented variants of the GVRSP.

\begin{sidewaystable}
\centering
\footnotesize
\caption{Some variants of the Green Vehicle Routing and Scheduling Problem} 
\label{TABLE3}

\begin{tabular}{|l|l|l|l|l|l|}

\\
\hline

Problem & $\mathcal{F}$ & $\mathcal{C}$ & Network characteristics & \multicolumn{2}{|c||}{Models and optimization methods} \\
 &  &  &  &  MultiObjective& SingleObjective \\
\hline
\\
RHM     & {Min. risk}    &(1) Origin-destination paths & Static network & Weighting method \cite{ABK88,CUR88,ROB81,SHO81,ZOG04} & Shortest path algorithms \cite{AKG02, JOH92, KUB97, LOM93}\\
        & Min. travel time                     & (2) {Risk equity} & {Stochastic dynamic} & Stochastic dominance \cite{TUR93} & Lagrangian relaxation \cite{GOP90a}\\

        & {Max. risk equity}            & (3) Vehicle capacity			     &{network}&Branch-and-bound \cite{KAR95}&Multi-commodity flow model \cite{CAR08}\\

	& Min. travel cost  & 			     &&Goal Programming \cite{GIA98}& Bilevel flow model \cite{BIA09}\\
	&           & 			     &&Insertion heuristic \cite{ZOG04}&\\
\hline
\\
TDVRP   & Min. travel time  &(1) Origin-destination paths & {The travel time or the}&Dynamic programming \cite{KOS93}& Nearest-neighbour heuristic \& \\

         & Min. travel cost &(2) Each costumer must be    & {travel speed on links}  & & branch-and-cut algorithm \cite{MAL92}\\

	 &					& serviced by one vehicle     &  {is a function of} &&Taboo search heuristic \cite{ICH03} \\
	 &					&(3) Vehicle capacity	      & {the time-of-day} &&Ant colony optimization \cite{DON08, VAN07}\\
&					&	      &  &&Local search algorithm \cite{HAS08}\\

\hline
\\
DVRP   & Min. travel time   &(1) Origin-destination paths & {The travel time on} & Ant colony optimization \cite{JUN08}&Tabu search algorithms \cite{VAN08}\\

        & Min. travel cost  &(2) Each costumer must be    & {links are stochastic}& Hybrid dynamic programming &Genetic algorithms \cite{TAN04} \\
	& Min. vehicle number	& serviced by one vehicle     & {varying over the time} & - ant colony optimisation \cite{CHI04}&Local search heuristics \cite{POT06}\\
	&					&(3) Vehicle capacity	      & & &\\

\hline
\\
RTVRP   & Min. travel time &(1) Origin-destination paths & {The travel time on} &&Genetic algorithm \cite{OKH09} \\
        &                  &(2) Each costumer must be    & {links are real}& &\\
	&					& serviced by one vehicle     & {time travel time} & &\\
	&					&(3) Vehicle capacity	      & &&\\

\hline
\\
CARP    & Min. travel time	&(1) Origin-destination paths & Static network & Genetic algorithm \cite{LAC03}& Tabu search algorithm \cite{GRE03, HER00}\\
	& Min. makespan 	&(2) Each link must be        & & Memetic algorithm \cite{MEI10}& Genetic and memetic algorithms \cite{LAC04}
\\
        & Min. travel cost 	& serviced by one vehicle     & & Epsilon-constraint method \cite{GRA10}& Guided local search algorithm \cite{BEU03}\\
	&			&(3) Vehicle capacity	      & && Ant colony optimization \cite{LAC04a}
\\

\hline
\\
MMVRP   & Min. travel time		& (1) Origin-destination paths  & Static multi-modal &Chance-constrained &  Double-sweep method \cite{BOR97}\\
	& {Max. cleaner}	&  (2) Vehicle capacity &  network        &goal programming \cite{MIN91}&Column generation \cite{MOC08}\\
	& {transport mean}		&   & &Subgradient optimization \cite{CHA08}&\\
	& Min. travel cost		&   & &  &\\

\hline
\\
DARP    & Max. costumers served	&(1) Origin-destination paths & Static network&Variable neighborhood-based heuristic \cite{PAR09}&Dynamic programming \cite{DES86, PSA80}\\
	& Min. vehicles number	&(2) A new costumer is accepted or & {Dynamic network} &&Branch-and-cut algorithm \cite{COR06a}\\
	& Max. level of service  	&rejected (3) Vehicle capacity  & &&\\
	&				&(4) Pairing and ride time	      & &&\\
\hline
\\

PDVRP    & Min. travel cost	&(1) Origin-destination paths  &  Static network &Insertion based algorithm \cite{MAD95}&Genetic algorithm \cite{JUN00}\\

    & Min. travel time	&(2) Capacity constraints  &  &&\\

    & Min. vehicles number	&(3) Precedence constraints  &  &&\\

\hline
\\

ERP	& {Min. energy used}		&(1) Origin-destination paths          & {The energy consumpt-} && Shortest-path algorithm \cite{ART10}\\
	&					&(2) {The battery cannot be}    & {ion/recuperation is} &&\\
	&					& {discharged below zero}       & {evaluted on each link}  & &\\
	&					&(3) {Battery energy capacity}  &  & &\\

\hline
\\
ATC	& Min. {delays}&(1) Origin-destination paths & Multi-dimentional network&Genetic algorithm \cite{TIA10}&Mixed-integer linear program- \\
	& Min. {fuel costs}	&(2) {Separation constraints} & &Hybrid ant colony optimization
 &optimization software \cite{PAL02}\\
	& Min. Max. {deviation}&(3) Bounds on travel time & &with a multi-agent approach \cite{VIJ08} & Lagrangian relaxation-heuristics \cite{BER00}\\
& &and distances of the new routes &  & \\
\hline
\end{tabular}
\end{sidewaystable}

Due to the conflicting nature of the criteria in sustainable transportation problems (economic, environmental and social), multiobjective optimization represents a major component. The particularity of these problems is that a unique feasible solution optimizing all the criteria does not exist. To obtain the optimal solutions, one just needs to consider Pareto optimal solutions. We can find in the literature two classes of multiobjective optimization methods, based on the problem elements type: \textit{deterministic methods} are the most studied in the literature and consider that all problem elements are deterministic and \textit{stochastic methods} consider that some elements of the problem are uncertain, these elements are modeled using random variables.

The unpredictable nature of transportation leads to many stochastic elements in the problem like the travel time, speed, traffic congestion, weather conditions, the amount of population present near a route and the effects of an accident. Deterministic methods use approximations of these elements and forecasts and sometimes lead to infeasible schedules and poor decisions. Recent research \cite{VAN06} shown the benefits of adaptability for vehicle routing and scheduling, exploiting stochastic information to produce better solutions. 

Generally, the uncertainty can be in the presence or absence of the customers, in the quantity of the their orders, and in the travel and service times. Routing and scheduling when demands are stochastic have been extensively studied when less attention is given to the case of stochastic travel time or stochastic speed. Especially with respect to travel times, variability is generally reduced to be nearly constant within time periods in a day. Such a characterization of travel times leads to the TDVRP variant (section \ref{SS_TDVRSP}). The stochastic version of multiobjective routing and scheduling is more studied in the context of transportation of hazardous materials \cite{CHA05, LEV92, MIL98, TUR93, WIJ93}, and essentially for the computation of shortest path problems. Very scarce work in the literature interest on the stochastic version of multiobjective vehicle routing and scheduling problem \cite{KEN03, LAP92}.  

As we can observe in table \ref{TABLE3}, heuristic and metaheuristic methods are the most used for solving multiobjective problems, we observe also that multiobjective problems are less studied than single objective problems (as our knowledge, no work in the literature leads with the multiobjective RTVRP and the multiobjective ERP). Evolutionary algorithms are one of the most popular methods for solving multiobjective routing problems \cite{JOZ08, JOZ02, TAN06}, these methods have been hybridized with local searches, heuristics, and/or exact methods for the problem resolution \cite{JOZ02, TAN06}. Other optimization methods were proposed in the literature for the resolution of multiobjective problems, based on genetic algorithms, lexicographic strategies, ant colony mechanisms, or specific heuristics \cite{JOZ08}, but a very limited works deal with the stochastic version of the problem.

\section{Conclusion}
\label{C_C}

This paper aims at specifying the contribution of Combinatorial Optimization (CO) to environmental transportation. To this purpose a framework for identifying relevant related problems treated in the CO field is proposed, some of the literature has been reviewed and discussed with respect to both models and optimization aspects. It can be observed that during the last few years, CO for transportation problems has extended its scope to include environmental applications.

It can be deduced that relationship between CO and environmental transportation is interactive in the sense that from the complexity of the issues examined stems the need to develop and adapt specific methodological tools. We synthesize in this section the material we have reviewed for the green routing and scheduling by summarizing some fundamental characteristics of this class of problems. 

\textit{A recent research area.} Green routing is a relatively recent problem since the consideration of environmental impacts in the models is a new issue. Several small research communities in this field work on their own problems, this causes a lack of common problem definitions, hypothesis, definitions and concepts. Decision-making in environmental transportation can be complex and seemingly intractable, principally because of the inherent trade-offs between socio-political, environmental, ecological, and economic factors. A balance has to be found between the complexity of the real world operation and the level at which the model is developed and also takes into account model accuracy and computational efficiency. 

\textit{Evaluation of environmental impacts.} Negative environmental impacts due to the transportation activity are multiple, we cite for example the land use for transport infrastructure, solid and hazardous waste and noise. The integration of the environmental costs of transportation is faced to many difficulties and is rarely quoted in the literature:

\begin{itemize}
 \item Fuel consumption and emissions are complex to estimate and are a function of several variables such as the type of vehicle, the speed, the acceleration rates and the meteorological conditions.
 \item Direct and indirect pollution have to be taken into account simultaneously for a good estimation of impacts. Direct pollution arises from the vehicles use (fuel) and indirect pollution arises from the energy generation such as electricity for example. 
 \item Pollution generated by transportation activity, once thought of as a purely local issue, now is recognized as a complex problem that is also subject to regional and global influences. 
\end{itemize}

These aspects of the problem have received very little interest in OR community, even if these are probably the most important issues surrounding transportation in the next century.

\textit{Multi-criteria approaches.} The green vehicle routing problem is a typical multiobjective problem. When we deal with multiple routes (which is the case in practice), multiobjective shortest path problems are of limited efficiency, because we relax the interaction between vehicles (whereas it uses the same network). More adapted models are global routing where a fleet of vehicles (identical or heterogeneous) have to be routed through the network. With environmental issues assuming greater importance, it is desirable to consider global multiobjective methods with uncertain and time-varying criterion such as travel time and speed.   

\textit{Uncertainty.} Many authors recognized the uncertain nature of some characteristics in transportation area (the number of incidents on a road and the travel time for example), these characteristics can be modeled by means of random variables whose distributions may vary over time. However, variable distributions are hard to determine in many cases due to the scarcity of data. The exact probabilistic expressions are usually too complicated, which results in the use of approximations for optimization. Hence, expertise is needed for understanding well probabilistic modeling to capture the important aspects of the activity, in addition, competence is also needed on optimization techniques to decide which approximations are necessary and which tools to use. Due to uncertainty, the path attribute values may not be simply additive across arcs in the path, causing the criterion to be non-order-preserving. This prevents the use of traditional dynamic programming techniques in the solution method. To deal with this, methods have to be developed that work with non-order-preserving criterion.


\begin{thebibliography}{00}


\bibitem{ABK88}
M. Abkowitz and P.D. Cheng, 
Developing a risk-cost framework for routing truck movements of hazardous materials,
Accident Analysis \& Prevention, 20(1):39-51, 1988.

\bibitem{ABK92}
M. Abkowitz, M. Lepofsky and P. Cheng, 
Selecting criteria for designating hazardous materials highway routes, 
Transportation Research Record 1333, 30-35, 1992.

\bibitem{AKG02}
V. Akgun, E. Erkut and R. Batta, 
On finding dissimilar paths, 
European Journal of Operational Research 121(2):232-246, 2000.

\bibitem{ALE00}
S.M. Alexander and N.M. Waters,
The effects of highway transportation corridors on wildlife: a case study of Banff National Park,
Transportation Research Part C: Emerging Technologies, 8(1-6):307-320, 2000.


\bibitem{ANC05}
D. Anciaux, D. Roy and S. Mirdamadi,
Un Mod\`ele de simulation d'une cha\^ine logistique avec prise en compte des moyens de transports intermodaux,
Proceedings of 5\`eme Conférence Internationale sur la Conception et Production Int\'egr\'ees, Casablanca, Maroc, 2005.

\bibitem{ANG02}
E. Angelelli and M.G. Speranza,
The periodic vehicle routing problem with intermediate facilities,
European Journal of Operational Research, 137:233-247, 2002.

\bibitem{ART10}
A. Artmeier, J. Haselmayr, M. Leucker and M. Sachenbacher,
The optimal routing problem in the context of battery-powered electric vehicles,
Workshop: CROCS at CPAIOR-10, Second International Workshop on Constraint Reasoning and Optimization for Computational Sustainability, Bologna, Italy, 2010.

%


\bibitem{BAL06}
R. Baldacci, L. Bodin and A. Mingozzi,
The multiple disposal facilities and multiple inventory locations rollon–rolloff vehicle routing problem,
Computers \& Operations Research 33:2667-2702, 2006.
         

\bibitem{BAT88}
R. Batta and S. Chiu,
Optimal obnoxious paths on a network: transportation of hazardous materials,
Operations Research, 36(1):84-92, 1988.


\bibitem{BER00}
D. Bertsimas and S.S. Patterson,
The traffic flow management rerouting problem in air traffic control: a dynamic network flow approach,
Transportation Science, 34(3):239-255, 2000. 

\bibitem{BEU03}
P. Beullens, L. Muyldermans, D. Cattrysse, and D.V. Oudheusden, 
A guided local search heuristic for the capacitated arc routing problem,
European Journal of Operational Research 147(3):629-643, 2003.

\bibitem{BIA09}
L. Biancoa, M. Caramia, and S. Giordani,
A bilevel flow model for hazmat transportation network design,
Transportation Research Part C: Emerging Technologies, 17(2):175-196, 2009.

\bibitem{BLO95}
J.M. Bloemhof-Ruwaard, P. Van Beek, L. Hordijk and L.N. Van Wassenhove,
Interactions between operational research and environmental management,
European journal of operational research, 85(2):229-243, 1995. 


\bibitem{BOD00}
L. Bodin, A. Mingozzi, R. Baldacci, and M. Ball,
The rollon-rollof vehicle routing problem,
Transportation Science, 34, 271-288, 2000.



\bibitem{BON04}
Y.M. Bontekoning, C. Macharis and J.J Trip,
Is a new applied transportation research field emerging? A review of intermodal rail-truck freight transport literature,
Transportation research. Part A, Policy and practice,  38(1):1-34, 2004.

\bibitem{BOR97}
B.S. Boardman, E.M. Malstrom, D.P. Butler and M.H Cole,
Computer assisted routing of intermodal shipments,
Computers \& Industrial Engineering, 33(1-2):311-314, 1997.


\bibitem{BOW87}
G.G. Brown, C.J. Ellis, G.W. Graves and D. Ronen,
Wide area dispatching of mobil tank trucks,
Interfaces 17, 107-120, 1987.


\bibitem{CAR08}
M. Caramia and P. Dell'Olmo, 
Multiobjective management in freight logistics: increasing capacity, service level and safety with optimization algorithms,
Springer London Ltd, 2008.

%

\bibitem{CHA97}
I. Chabini, 
A New shortest paths algorithm for discrete dynamic networks,
Proceeding of 8th IFAC Symposium on Transport Systems, 551-556, 1997.


\bibitem{CHA08}
T.-S. Chang,  
Best routes selection in international intermodal networks,
Computers and Operations Research, 35(9):2877-2891, 2008. 

\bibitem{CHA05}
T. Chang, L. Nozick and M. Turnquist, 
Multiobjective path finding in stochastic dynamic networks, with application to routing hazardous materials shipments, 
Transportation Science 39(3):383-399, 2005.

\bibitem{CHI04}
D.M. Chitty and M.L. Hernandez,
A hybrid ant colony optimisation technique for dynamic vehicle routing,
Lecture Notes in Computer Science, 3102:48-59, 2004.

\bibitem{IRN05}
S. Irnich and G. Desaulniers,
Shortest Path Problems with Resource Constraints,
In G. Desaulniers, J. Desrosiers, and M.M. Solomon, editors, Column Generation, chapter 2, pages 33-65, Springer 2005. 

\bibitem{CLA08} 
J-P. Clarke, M. Lowther, L. Ren, W. Singhose, S. Solak, A. Vela and L. Wong,
En route traffic optimization to reduce environmental impact,
PARTNER Project 5 report, 2008.

\bibitem{COR06a}
J.-F. Cordeau, 
A branch-and-cut algorithm for the dial-a-ride problem, 
Operations Research, 54(3):573-586, 2006. 


\bibitem{COR06}
J.-F. Cordeau and G. Laporte,
The dial-a-ride problem (DARP): variants, modeling issues and algorithms,
4OR, 1:89-101, 2003.


\bibitem{CUR88}
J. R. Current and C. ReVelle and J. L Cohon,
The minimum-covering/shortest-path problem,
Decision Sciences, 19(3):490-503, 1988.


\bibitem{DAN97}
S.E. Daniel, D.C. Diakoulaki and C.P. Pappis,
Operations research and environmental planning,      
European journal of operational research, 102:248-263, 1997. 

\bibitem{DEM97}
L. De Meulemeester, G. Laporte, FV. Louveaux and F. Semet,
Optimal sequencing of skip collections and deliveries,
Journal of the Operational Research Society, 48:57-64, 1997


\bibitem{DES86}
J. Desrosiers, Y. Dumas and F. Soumis,
A dynamic programming solution of the largescale single-vehicle dial-a-ride problem with time windows,
American Journal of Mathematical and Management Sciences, 6:301-325, 1986.

\bibitem{DON08}
A.V. Donati, R. Montemanni, N. Casagrande, A.E. Rizzoli and L.M. Gambardella,
Time dependent vehicle routing problem with a multi ant colony system,
European Journal of Operational Research, 185(3):1174-1191, 2008.

\bibitem{DRO00}
M. Dror, 
Arc routing, theory, solutions, and applications. 
Boston, MA: Kluwer, 2000. 


\bibitem{EC96}
European Commission, 
State of the art on strategic environmental assessment for transport infrastructure,
Report by the European Commission, Brussels, 1996.

%

\bibitem{ERK05}
E. Erkut and A. Ingolfsson,
Transport risk models for hazardous materials: revisited,
Operations Research Letters, 33:81-89, 2005.

\bibitem{ERK98}
E. Erkut and V. Verter, 
Modeling of transport risk for hazardous materials, 
Operations Research, 46(5):625-642, 1998.



\bibitem{FLE97}
M. Fleischmann, J.M. Bloemhof-Ruwaard, R. Dekker, E. van der Laan, J.A.E.E. van Nunen, and L.N. van Wassenhove,
Quantitative models for reverse logistics: A review,
European journal of operational research, 103:1-17, 1997.

\bibitem{FLE05}
G. Fleury, P. Lacomme, C. Prins and W. Ramdane-Cherif, 
Improving robustness of solutions to arc routing problems,
Journal of the operational research society, 56:526-538, 2005.




%

\bibitem{GIA03}
G. Ghiani, F. Guerriero, G. Laporte and R. Musmanno,
Real-time vehicle routing: solution concepts, algorithms and parallel computing strategies,
European Journal of Operational Research, 151(1):1-11, 2003.

\bibitem{GIA98}
I. Giannikos,
A multiobjective programming model for locating treatment sites and routing hazardous wastes,
European Journal of Operational Research, 104:333-342, 1998.



\bibitem{GOP90a}
R. Gopalan, R. Batta and M.H. Karwan, 
The equity constrained shortest path problem, 
Computers and Operations Research, 17:297-307, 1990.

\bibitem{GOP90b}
R. Gopalan, K.S. Kolluri, R. Batta and M.H. Karwan, 
Modeling equity of risk in the transportation of hazardous materials, 
Operations Research, 38(6):961-973, 1990.

\bibitem{GRA10}
L. Grandinetti, F. Guerriero, D. Lagan\`a and O. Pisacane,
An approximate epsilon-constraint method for the multiobjective undirected capacitated arc routing problem, 
Lecture Notes in Computer Science, 214-225, 2010.

\bibitem{GRE03}
P. Greistorfer, 
A tabu scatter search metaheuristic for the arc routing problem,
Computers \& Industrial Engineering, 44(2):249-266, 2003.



%
%

\bibitem{HAR93}
D.W. Harwood, J.G. Viner and E.R. Russell,
Procedure for developing truck accident and release rates for hazmat routing,
J. Transp. Eng. 119:189-199, 1993.

\bibitem{HAS08}
H. Hashimoto, M. Yagiura and T. Ibaraki,
An iterated local search algorithm for the time-dependent vehicle routing problem with time windows,
Discrete Optimization, 5(2):434-456, 2008.

\bibitem{HIL92}
A.V. Hill and W.C. Benton,
Modeling intra-city time-dependent travel speeds for vehicle scheduling problems,
Journal of the Opertions Research Society, 43(4):343-351, 1992.

\bibitem{HER00}
A. Hertz, G. Laporte and M. Mittaz, 
A tabu search heuristic for the capacitated arc routing problem,
Operations Research, 48(1):129-135, 2000.

\bibitem{INS91}
S. Brian,  Bochner, Chairperson, Washington,
D.C.
Institute of Transportation Engineers,
Traffic access and impact studies for site development, 
A Recommended Practice, Transportation Planners Council Task Force on Traffic Access/Impact Studies, 

\bibitem{IYE09}
M.H. Iyengar, 
Incorporating an economic model in the health impact assessment approach. In: The 10th Annual Health Impact Assessment Conference, Rotterdam, Netherlands, 2009.


\bibitem{ICH03}
S. Ichoua, M. Gendreau and J-Y. Potvin,
Vehicle dispatching with time-dependent travel times,
European Journal of Operational Research, 144(2):379-396, 2003.  

\bibitem{JOH92}
P.E. Johnson, D.S. Joy, D.B. Clarke and J.M. Jacobi, 
HIGWAY 3.01, An enhanced highway routing model: Program, description, methodology and revised user's manual, Oak Ridge National Laboratory, ORLN/TM-12124, Oak Ridge, TN, 1992.

\bibitem{JOZ08}
N. Jozefowiez, F. Semet and El-Ghazali Talbi,
Multiobjective vehicle routing problems,
European Journal of Operational Research, 189(2):293-309, 2008.

\bibitem{JOZ02}
N. Jozefowiez, F. Semet and El-Ghazali Talbi,
Parallel and hybrid models for multi-objective optimization: application to the vehicle routing problem,
in: J.J. Merelo Guervos et al. (Eds.), Parallel Problem Solving from Nature VII, Lecture Notes in Computer Science, 2439/2002:271-280, 2002.


\bibitem{JUN08}
Q. Jun, J. Wang, and B. Zheng, 
A hybrid multi-objective algorithm for dynamic vehicle routing problems.
Lecture Notes in Computer Science, 5103:674-681, 2008.

\bibitem{JUN00}
S. Jung and A. Haghani, 
Genetic algorithm for a pickup and delivery problem with time windows, 
Transportation research record, 1733:1-7, 2000.



\bibitem{KEN03}
A.S. Kenyon and D.P. Morton,
Stochastic vehicle routing with random travel times,
Transportation Science, 37(1):69-82, 2003.

\bibitem{KIM06}
B.I. Kim, S. Kim and S. Sahoo,
Waste Collection Vehicle Routing Problem with Time Windows,
Computers \& Operations Research, 33:3624-2642, 2006.


\bibitem{KOS93}
MM. Kostreva and MM. Wiecek,
Time dependency in multiple objective dynamic programming,
Journal of Mathematical Analysis and Applications 173(1):289-307, 1993.


\bibitem{KOU93}
H.N. Koutsopoulos and H. Xu,
An information discounting routing strategy for advanced traveler information systems. 
Transportation Research Part C 1 (3):249-264, 1993.

\bibitem{KUB97}
M. Kuby, X. Zhongyi and X. Xiaodong, 
A minimax method for finding the $k$ best differentiated paths, 
Geographical Analysis 29(4):298-313, 1997.

\bibitem{KUC00}
J. Kuchar and L. Yang, 
A Review of conflict detection and resolution modeling methods,
IEEE Transactions on Intelligent Transportation Systems, vol. 1, 2000.


\bibitem{LAC03}
P. Lacomme, C. Prins and M. Sevaux, 
Multiobjective capacitated arc routing problem,
LNCS 2632, pp. 550-564, Springer, 2003.

\bibitem{LAC04a}
P. Lacomme, C. Prins and A. Tanguy,
First competitive ant colony scheme for the CARP,
Lecture Notes in Computer Science, Publisher: Springer-Verlag Heidelberg, 3005:426-427, 2004.


\bibitem{LAC04}
P. Lacomme, C. Prins and W. Ramdane-Cherif, 
Competitive memetic algorithms for arc routing problems,
Annals of Operations Research, 131:159-185, 2004.

\bibitem{LAP92}
G. Laporte, F. Louveaux and H. Mercure,
The vehicle routing problem with stochastic travel times,
Transportation Science, 26(3):161-170, 1992.

\bibitem{LEV92}
H. Levy,
Stochastic dominance and expected utility: Survey and analysis. 
Management Science 38(4):555-593, 1992.

\bibitem{LIS91}
G.F. List, P.B. Mirchandani, M.A. Turnquist and K.G. Zografos, 
Modeling and analysis for hazardous materials transportation: Risk analysis, routing/scheduling and facility location, 
Transportation Science, 25(2):100-114, 1991.


\bibitem{LIT06}
T. Litman,
Transportation Market Distortions,
Berkeley Planning Journal, 19:19-36, 2006.

\bibitem{LOM93}
K. Lombard and R.L. Church, 
The gateway shortest path problem: Generating alternative routes for a corridor location problem, 
Geographical Systems, 1:25-45, 1993.


\bibitem{LEC09}
C. Lecluyse, T. Van Woensel and H. Peremans,
Vehicle routing with stochastic time-dependent travel times,
4OR, 7(4): 363-377, 2009.


\bibitem{MAD95}
O.B.G. Madsen, H.F. Ravn and J.M. Rygaard,
A heuristic algorithm for a dial-a-ride problem with time windows, multiple capacities, and multiple objectives, Annals of Operations Research, 60:193-208, 1995.

\bibitem{MAL92}
C. Malandraki and M.S. Daskin,
Time dependent vehicle routing problems: formulations, properties and heuristic algorithms,
Transportation Science, 26:185-200, 1992.


\bibitem{MEI10}
Y. Mei, K. Tang and X. Yao, 
Decomposition-based memetic algorithm for multiobjective capacitated arc routing problem,
IEEE Transactions on Evolutionary Computation, to appear, 2010. 

\bibitem{MIL01}
E.D. Miller-Hooks, 
Adaptive least-expected time paths in stochastic, time-varying transportation and data networks,
Networks, 37(1):35-52, 2001.

\bibitem{MIL98}
E.D. Miller-Hooks and H.S. Mahamassani,  
Optimal routing of hazardous materials in stochastic, time-varying transportation networks,
Transportation Research Record, 1645, 143-151, 1998.

\bibitem{MIN91}
H. Min,
International intermodal choices via chance-constrained goal programming,
Transportation Research Part A: General, 25(6):351-362, 1991.


\bibitem{MOC08}
L. Moccia, J. F. Cordeau, G. Laporte, S. Ropke and M.P. Valentini,
Modeling and solving a multimodal routing problem with timetables and time windows. 
submitted to Networks, 2008.



\bibitem{NEM99}
A.K. Nemaa and S.K. Gupta,
Optimization of regional hazardous waste management systems: an improved formulation,
Waste Management, 19:441-451, 1999.

\bibitem{NGU86}
S. Nguyen and S. Pallottino,
Hyperpaths and shortest hyperpaths,
In: Combinatorial Optimization. Lecture Notes in Mathematics, 1403:258-271, Berlin, 1986.

\bibitem{NUO06}
T. Nuortio, J. Kytojokib, H. Niskaa and O. Braysy,
Improved route planning and scheduling of waste collection and transport,
Expert Systems with Applications, 30:223-232, 2006.



\bibitem{OFF96}
Office of Policy, Planning and Evaluation,
Indicators of the Environmental Impacts of Transportation: Highway, Rail, Aviation, and Maritime Transport,
U.S. Environmental Protection Agency, Washington, D.C, 1996.

\bibitem{OKH09}
I. Okhrin and K. Richter,
Vehicle routing problem with real-time travel times,
International journal of vehicle information and communication systems, 2(1-2):59-77, 2009.



\bibitem{PAR09}
S.N. Parragh, K.F. Doerner, R.F. Hartl and X. Gandibleux, 
A heuristic two-phase solution approach for the multi-objective dial-a-ride problem,
Networks 54(4):227-242, 2009.

\bibitem{PAR08}
S.N. Parragh, K.F. Doerner and R.F. Hartl,
A survey on pickup and delivery problems,
part i: Transportation between customers and depot, Journal f\"{u}r Betriebswirtschaft, 58(1):21-51, 2008.
          
\bibitem{PAL02}
L. Pallottino, E. M. Feron and A. Bicchi, 
Conflict resolution problems for air traffic management systems solved with mixed integer programming,
IEEE Transactions on Intelligent Transportation Systems, vol. 3, 2002.

\bibitem{PAL07}
A. Palmer, 
The development of an integrated routing and carbon dioxide emissions model for goods vehicles, 
Cranfield, PhD thesis, 2007. 

\bibitem{POT06}
J-Y. Potvin, Y. Xua and I. Benyahia,
Vehicle routing and scheduling with dynamic travel times,
Computers and Operations Research, 33(4):1129-1137, 2006.

\bibitem{PSA95}
H.N. Psaraftis,
Dynamic vehicle routing: status and prospects,
Annals of Operations Research, 61:143-164, 1995.

\bibitem{PSA80}
H.N. Psaraftis,
A dynamic programming solution to the single vehicle many-to-many immediate request dial-a-ride problem,
Transportation Science, 14(2):130-154, 1980.


\bibitem{RAV91}
C. ReVelle, J. Cohon and D. Shobrys,
Simultaneous Siting and Routing in the Disposal of Hazardous Wastes,
Transportation Science, 25(2):138-145, 1991.

\bibitem{ROB81}
J.C. Robbins, 
Routing hazardous materials shipments,
Thesis/Dissertation, Indiana Univ.,Bloomington, IN, 1981.


\bibitem{SAC85}
F.F Saccomanno and A. Chan,
Economic evaluation of routing strategies for hazardous road shipments,
Transportation Research, 1020:12-18, 1985.

\bibitem{SAH05}
S. Sahoo, S. Kim, B-I. Kim, B. Kraas and A. Popov, 
Routing Optimization for Waste Management,
Interfaces, 35(1):24-36, 2005.

\bibitem{SBI07}
A. Sbihi and R.W. Eglese,
Combinatorial optimization and green logistics,
4OR, 5:99-116, 2007.

\bibitem{SCH09}
D. Schram-Bijkerk, E. Van Kempen, A.B. Knol, H. Kruize, B. Staatsen and I. Van Kamp,
Quantitative health impact assessment of transport policies: two simulations related to speed limit reduction and traffic re-allocation in the Netherlands,
Occupational and environmental medicine, 66(10):691-698, 2009.

\bibitem{SHO81}
D. Shobrys, 
A model for the selection of shipping routes and storage locations for a hazardous substance,
Ph.D Dissertation, John Hopkins University, Baltimore, 1981.

\bibitem{SIV95}
R.A. Sivakumar, R. Batta and M.H. Karwan, 
A multiple route conditional risk model for transporting hazardous materials, 
INFOR 33:20-33, 1995.


%




\bibitem{TAN04}
E. Taniguchi and H. Shimamoto,
Intelligent transportation system based dynamic vehicle routing and scheduling with variable travel times,
Transportation Research Part C: Emerging Technologies, 12(3-4):235-250, 2004.

\bibitem{TAN06}
K.C. Tan, Y.H. Chew, L.H. Lee,
A hybrid multi-objective evolutionary algorithm for solving truck and trailer vehicle routing problems, 
European Journal of Operational Research, 172:855-885, 2006.

\bibitem{TAR07}
Z. Tarapata, 
Selected multicriteria shortest path problems : An analysis of complexity, models and adaptation of standard algorithms,
International Journal Of Applied Mathematics And Computer Science, 17:269-287, 2007.


\bibitem{TIA10}
W. Tian and M. Hu, 
Study of air traffic flow management optimization model and algorithm based on multi-objective programming,
Second International Conference on Computer Modeling and Simulation, 2:210-214, 2010.

\bibitem{TOTB02}
P. Toth and D. Vigo,
The vehicle routing problem, 
Philadelphia: Society for Industrial and Applied Mathematics, 2002.

\bibitem{TSAB00}
D. Tsamboulas and G. Mikroudis,
EFECT - evaluation framework of environmental impacts and costs of transport initiatives, 
Transportation Research Part D (5):283-303, 2000.

\bibitem{TUN00}
DV. Tung and A. Pinnoi,
Vehicle routing–scheduling for waste collection in Hanoi,
European Journal of Operational Research, 125(3):449-468, 2000. 

\bibitem{TUR93}
M. Turnquist, 
Multiple objectives, uncertainty and routing decisions for hazardous materials shipments, 
Proceedings of the 5th International Conference on Computing in Civil and Building Engineering ASCE, New York 357-364, 1993.




\bibitem{VAN06}
P. Van Hentenryck and R. Bent,
Online Stochastic Combinatorial Optimization,
MIT Press, 2006.

\bibitem{VAN07}
T. Van Woensel, L. Kerbache, H. Peremans and N. Vandaele,
A queueing framework for routing problems with time-dependent travel times,
Journal of Mathematical Models and Algorithms, Special Issue ``Quantitative aspects of transportation and logistics'', 6(1):151-173, 2007.

\bibitem{VAN08}
T. Van Woensel and L. Kerbache and H. Peremans and N. Vandaele,
Vehicle routing with dynamic travel times : A queueing approach,
European journal of operational research, 186(3): 990-1007, 2008. 

\bibitem{VEE05}
J.L. Veerman, J.J. Barendregt and J.P. Mackenbach,
Quantitative health impact assessment: current practice and future directions,
Journal of Epidemiology and Community Health, 59(5):361-370, 2005. 

\bibitem{VIJ08}
V.P. Vijayan, D. John, M. Thomas, N.V. Maliackal and S.S. Vargheese,
Multi agent path planning approach to dynamic free flight environment,           
International Journal of Recent Trends in Engineering (IJRTE), 1(1):41-46, 2009.




\bibitem{WAN10}
Y.W Wang and C.R. Wang,
Locating passenger vehicle refueling stations,
Transportation Research Part E: Logistics and Transportation Review, 46(5):791-801, 2010. 

\bibitem{WIJ93}
A.B. Wijeratne, M.A. Turnquist and P.B. Mirchandani, 
Multiobjective routing of hazardous materials in stochastic networks, 
European Journal of Operational Research, 65:33-43, 1993.





\bibitem{KAR95}
J. Karkazis and T.B. Boffey,
Optimal location of routes for vehicles transporting hazardous materials,
European journal of operational research, 86(2):201-215, 1995.	

\bibitem{ZHU00}
K.Q. Zhu and K.-L. Ong,
A reactive method for real time dynamic vehicle routing problems, 
In Proceedings of 12th IEEE Internationals Conference on Tools with Artificial Intelligence. IEEE Computing Society, Los Alamitos, USA, 176-180, 2000. 

\bibitem{ZOG04}
K.G. Zografos and K.N. Androutsopoulos,
A heuristic algorithm for solving hazardous materials distribution problems,
European Journal of Operational Research, 152(2):507-519, 2004. 

\bibitem{ZOG89}    
K.G Zografos and C.F. Davis, 
Multiobjective programming approach for routing hazardous materials,
Journal of Transportation Engeneering 115(6):661-673, 1989.





\end{thebibliography}
\end{document}